\documentclass[12pt]{article}
\usepackage{amssymb}
\usepackage{amsfonts}
\usepackage{amsmath}
\usepackage{multirow}
\usepackage{titlesec}

\titleformat*{\section}{\large\bfseries}
\titleformat*{\subsection}{\normalsize\bfseries}
\titleformat*{\subsubsection}{\normalsize\bfseries}
\titleformat*{\paragraph}{\normalsize\bfseries}
\titleformat*{\subparagraph}{\normalsize\bfseries}

\newtheorem{theorem}{Theorem}
\newtheorem{proposition}[theorem]{Proposition}
\newtheorem{lemma}[theorem]{Lemma}

\begin{document}

\title{Circles in self dual symmetric $R$-spaces}
\author{Marcos Salvai\thanks{
This work was partially supported by Conicet (PIP 112-2011-01-00670), Foncyt
(PICT 2010 cat 1 proyecto 1716) and Secyt Univ.\thinspace Nac.\thinspace C 
\'{o}rdoba.} \\
{\small {CIEM - FaMAF, Conicet - Universidad Nacional de C\'ordoba }}\\
{\small {Ciudad Universitaria, 5000 C\'{o}rdoba, Argentina}} \\
{\small {salvai@famaf.unc.edu.ar}}}
\date{ }
\maketitle

\begin{abstract}
Self dual symmetric $R$-spaces have special curves, called circles, 
introduced by Burstall, Donaldson, Pedit and Pinkall in 2011, whose 
definition does not involve the choice of any Riemannian metric. We 
characterize the elements of the big transformation group $G$ of a self dual 
symmetric $R$-space $M$ as those diffeomorphisms of $M$ sending circles in 
circles. Besides, although these curves belong to the realm of the 
invariants by $G$, we manage to describe them in Riemannian geometric 
terms:\ Given a circle $c$ in $M$, there is a maximal compact subgroup $K$ 
of $G$ such that $c$, except for a projective reparametrization, is a 
diametrical geodesic in $M$ (or equivalently, a diagonal geodesic in a 
maximal totally geodesic flat torus of $M$), provided that $M$ carries the 
canonical symmetric $K$-invariant metric. We include examples for the 
complex quadric and the split standard or isotropic Grassmannians. 
\end{abstract}


\medskip


\section{Introduction}

Symmetric $R$-spaces are compact Riemannian symmetric spaces admitting a Lie
group of diffeomorphisms, called the big transformation group, properly 
containing the isometry group. It is a celebrated theorem of Nagano \cite%
{Nagano} that, up to covers, these are essentially the only compact 
Riemannian symmetric spaces admitting such a group.

If the action of a group $G$ on a manifold $M$ is given, it is natural to 
look for a structure on $M$ whose automorphism group is $G$. For instance, 
the sphere $S^{n}$ admits two different big transformation groups, those of 
conformal and projective maps, respectively. The big transformation group of
a compact Hermitian symmetric space $M$ is the group of its 
(anti-)holomorphic maps. For a general symmetric $R$-space $M$, the big 
transformation group has been identified as the group of diffeomorphisms 
preserving a geometrically defined quantity on $M$ called the arithmetic 
distance \cite{BasicTakeuchi}, or preserving a so called generalized 
conformal structure (the smooth assignment of a cone in each tangent space 
of $M$, satisfying certain properties) \cite{GK} (see also \cite{Bertram}).

Self dual symmetric $R$-spaces have special curves, called circles, 
introduced in Burstall et al \cite{BP} in order to extend some aspects of 
the theory of Darboux transforms of isothermic surfaces in $\mathbb{R}^{3}$ 
to significantly more general ambient spaces. Their definition does not 
involve the choice of any Riemannian metric. In the case of the sphere $
S^{n} $, their trajectories are the usual circles (intersections of $S^{n}$ 
with planes at distance smaller than one from the origin). See the list of 
all self dual symmetric spaces at the end of \cite{BP} and detailed examples
in Subsection 2.2 of that article and Section \ref{Examples} below.

In Theorem \ref{Uno} we characterize the elements of the big transformation
group of a self dual symmetric $R$-space $M$ as those diffeomorphisms of $M$
sending circles in circles. Theorem \ref{Dos} gives some geometrical
properties of circles. Although circles in a self dual symmetric $R$-space $M
$ belong to the realm of the invariants of the big transformation group of $M
$, we manage to describe them in Riemannian geometric terms:\ Given a circle 
$c$ in $M$, there is a maximal compact subgroup $K$ of $G$ such that, except
for a projective reparametrization, $c$ is a diametrical geodesic in $M$ (or
equivalently, a diagonal geodesic in a maximal torus of $M$), provided that $%
M$ carries the canonical symmetric $K$-invariant metric.

In \cite{LO'H}, Langevin and O'Hara presented a new approach to the
classical subject of conformal length of curves in the sphere, in terms of
the canonical pseudo Riemannian structure on the space $\mathcal{C}$ of
oriented circles: The circles osculating a curve $\alpha $ in the sphere
define a null curve in $\mathcal{C}$, whose $\frac{1}{2 }$-dimensional length
provides, generically, a conformally invariant parametrization of $\alpha $.
The curves we deal with in this article seem to be adequate to generalize
that result to self dual symmetric $R$-spaces. We are currently working on
the subject.

We adopt the approach of \cite{BP}, only taking for granted the action of 
the big transformation group of $M$. We emphasize that no Riemannian metric 
on $M$ is fixed from the beginning, in contrast to \cite{BasicTakeuchi, GK}.
Since our proofs are based on their results, compatible symmetric Riemannian
metrics are necessarily auxiliary constituents, even if they are not 
involved in some statements.

Next we recall from \cite{BP} (Sections 1 and 2 and Subsection 4.1) the 
definitions and basic properties of symmetric $R$-spaces, self duality and 
circles.

\subsection{Self dual symmetric $R$-spaces}

Let $\mathfrak{g}$ be a real simple Lie algebra. Let Inn~$\left( \mathfrak{g}
\right) $ be the unique connected Lie subgroup of $GL\left( \mathfrak{g}
\right) $ whose Lie algebra is ad~$\left( \mathfrak{g}\right) =\left\{ \text{
ad}_{x}\mid x\in \mathfrak{g}\right\} \subset $ gl$~\left( \mathfrak{g}
\right) $. For $g\in $ Inn~$\left( \mathfrak{g}\right) $, let $I_{g}$ be the
automorphism of Inn~$\left( \mathfrak{g}\right) $ defined as usual by $
I_{g}\left( h\right) =ghg^{-1}$, satisfying Ad$~\left( g\right) =\left( 
dI_{g}\right) _{e}:$ ad~$\left( \mathfrak{g}\right) \rightarrow $ ad~$\left(
\mathfrak{g}\right) $. We may think of Inn~$\left( \mathfrak{g}\right) $ as 
the adjoint group of $\mathfrak{g}$, since ad $:\mathfrak{g}\rightarrow $
ad~ $\left( \mathfrak{g}\right) $ is an isomorphism ($\mathfrak{g}$ is
simple).

\bigskip

\noindent \textbf{Definitions. }Let\textbf{\ }$\mathfrak{g}$ be a simple 
real Lie algebra. A subalgebra $\mathfrak{p}$ of $\mathfrak{g}$ is said to 
be \emph{height one parabolic} if the polar $\mathfrak{p}^{\bot }=\left\{ 
x\in \mathfrak{g}\mid B\left( x,y\right) =0\text{ for all }y\in \mathfrak{p}
\right\} $ with respect to the Killing form $B$ of $\mathfrak{g}$ is an 
abelian subalgebra of $\mathfrak{g}$. Notice that $\mathfrak{p}^{\bot
}\subset \mathfrak{p}$.

A \emph{symmetric }$R$\emph{-space} is a conjugacy class of height one 
parabolic subalgebras of $\mathfrak{g}$, that is, the orbit of a height one 
parabolic subalgebra of $\mathfrak{g}$ under the action of Inn~$\left(  
\mathfrak{g}\right) $.

Two height one parabolic subalgebras $\mathfrak{p},\mathfrak{q}$ of $ 
\mathfrak{g}$ are \emph{opposite} if $\mathfrak{p}^{\bot }\cap \mathfrak{q}
^{\bot }=\left\{ 0\right\} $, or equivalently, if $\mathfrak{g}=\mathfrak{p}
\oplus \mathfrak{q}^{\bot }$. We adopt the nomenclature of \cite{Frost}; in  
\cite{BP} they are said to be complementary.

The \emph{dual} $M^{\ast }$ of a symmetric $R$-space $M$ is defined to be 
the set of parabolic subalgebras $\mathfrak{q}$ of $\mathfrak{g}$ for which $
\mathfrak{q}$ is opposite to some $\mathfrak{p}\in M$. It turns out to be a 
symmetric $R$-space as well. One says that $M$ is \emph{self-dual} if $
M^{\ast }=M$, that is, if for any $\mathfrak{p}\in M$, any height one 
parabolic subalgebra $\mathfrak{q}$ of $\mathfrak{g}$ opposite to $\mathfrak{%
\ p}$ is in $M$. Equivalently, if for some $\mathfrak{p}\in M$, some height 
one parabolic subalgebra $\mathfrak{q}$ of $\mathfrak{g}$ opposite to $ 
\mathfrak{p}$ is conjugate to $\mathfrak{p}$.

\subsection{Circles in self dual symmetric $R$-spaces\label{SScircles}}

Let $\mathfrak{q}$ be a height one parabolic subalgebra of a simple real Lie
algebra $\mathfrak{g}$. Then $\left. \exp \right\vert _{\mathfrak{q}^{\bot
}}:\mathfrak{q}^{\bot }\rightarrow \exp \left( \mathfrak{q}^{\bot }\right) $
is a global diffeomorphism onto the abelian group $\exp \left( \mathfrak{q}
^{\bot }\right) $. Let $M$ be a symmetric $R$-space and let $\mathfrak{q}\in
M$. The set  
\begin{equation*}
\Omega _{\mathfrak{q}}=\left\{ \mathfrak{p}\in M\mid \mathfrak{p}\text{ is
opposite to }\mathfrak{q}\right\}
\end{equation*}
is open and dense in $M$. Moreover, $\exp \left( \mathfrak{q}^{\bot }\right)
$ acts simply transitively on it. Hence, given $\mathfrak{p}\in \Omega _{ 
\mathfrak{q}}$, the map  
\begin{equation}
\mathfrak{q}^{\bot }\rightarrow \Omega _{\mathfrak{q}}\text{,\ \ \ \ \ \ \ \
\ \ \ \ \ }x\mapsto \left( \exp x\right) \cdot \mathfrak{p}  \label{stereo}
\end{equation}
is a diffeomorphism. It is called inverse-stereoprojection with respect to $
\left( \mathfrak{p},\mathfrak{q}\right) $. This induces a linear isomorphism
\begin{equation}
\pi _{\mathfrak{q}}^{\mathfrak{p}}:\mathfrak{q}^{\bot }\rightarrow T_{ 
\mathfrak{p}}M,\ \ \ \ \ \ \ \ \ \ \ \ \ \ \pi _{\mathfrak{q}}^{\mathfrak{p}
}\left( x\right) =\left. \frac{d}{dt}\right\vert _{0}\exp \left( tx\right)
\cdot \mathfrak{p}  \label{TpMqBot}
\end{equation}
for any $\mathfrak{p}\in \Omega _{\mathfrak{q}}$. This identification 
remained mostly implicit in \cite{GK}. We make it explicit for the sake of 
clarity of the exposition.

\bigskip

\noindent \textbf{Definition. } Let $M$ be a self dual symmetric $R$-space.
Let $\mathfrak{p },\mathfrak{p}_{1}$ and $\mathfrak{q}$ be three points in $M
$ which are  pairwise opposite. Since the map in (\ref{stereo}) is a
diffeomorphism, then  there exists a unique $y\in \mathfrak{q}^{\bot }$ such
that $\mathfrak{p} _{1}=\exp \left( y\right) \cdot \mathfrak{p}$. The curve $%
c:\mathbb{R}\cup  \left\{ \infty \right\} \rightarrow M$ 
\begin{equation*}
c\left( t\right) =\exp \left( ty\right) \cdot \mathfrak{p},\ \ \ \ \ \ \ \
c\left( \infty \right) =\mathfrak{q}
\end{equation*}
is called the \emph{circle in }$M$\emph{\ through} $\mathfrak{p},\mathfrak{p}
_{1},\mathfrak{q}$.

\medskip

The image of $c$ is in fact a circle. Moreover, if the circles determined by
two triples of (appropriately ordered) pairwise  opposite points have the
same trajectories, then their parametrizations  differ in a fractional
linear transformation. See the precise assertion in  Proposition 4.4 of \cite%
{BP}.

\subsection{The big transformation group}

Now we recall from \cite{GK} the definition of a Lie group $G$ acting on $M$
which is in  general bigger than Inn~$\left( \mathfrak{g}\right) $. For
instance, for the conformal structure of $ S^{n}$, $G$ will be the group $%
O_{\uparrow }\left( 1,n+1\right) $ of all  conformal maps, not only the
group Inn~$\left( o\left( 1,n+1\right) \right)  =O_{0}\left( 1,n+1\right) $
consisting of those conformal maps preserving  the orientation.

Let $M$ be a self dual symmetric $R$-space. Given two opposite elements $ 
\mathfrak{p},\mathfrak{q}$ in $M$, we consider the group  
\begin{equation*}
G_{0}\left( \mathfrak{p},\mathfrak{q}\right) =\left\{ g\in \text{Aut~}\left( 
\mathfrak{g}\right) \mid g\left( \mathfrak{p}\right) =\mathfrak{p}\text{ and 
}g\left( \mathfrak{q}\right) =\mathfrak{q}\right\}
\end{equation*}
and call  
\begin{equation*}
G=G_{0}\left( \mathfrak{p},\mathfrak{q}\right) ~\text{Inn}~\left( \mathfrak{%
g }\right) \text{,}
\end{equation*}
which is a subgroup of Aut~$\left( \mathfrak{g}\right) $ with the same Lie 
algebra ad$~\left( \mathfrak{g}\right) $ as Inn~$\left( \mathfrak{g}\right) $
. This follows from the fact that if $h\in G_{0}\left( \mathfrak{p}, 
\mathfrak{q}\right) $ and $g\in $ Inn~$\left( \mathfrak{g}\right) $, then $
hgh^{-1}\in $ Inn~$\left( \mathfrak{g}\right) $.

By Proposition 2.5 in \cite{BP}, the group Inn~$\left( \mathfrak{g}\right) $
acts transitively on the pairs of opposite elements of $M$. Hence, $
G_{0}\left( \mathfrak{p},\mathfrak{q}\right) $ and $G_{0}\left( \mathfrak{p}
^{\prime },\mathfrak{q}^{\prime }\right) $ are conjugate by an element of 
Inn~$\left( \mathfrak{g}\right) $. Therefore, the group $G$ depends only of 
the symmetric $R$-space $M$, and not on the choice of opposite $\mathfrak{p}%
, \mathfrak{q}$ in $M$. It is called the \emph{big transformation group} of $%
M$ .

The following theorem characterizes the big transformation group of a self 
dual symmetric $R$-space.

\begin{theorem}
\label{Uno}Let $M$ be a self dual symmetric $R$-space and let $G$ be its big
transformation group. Then a diffeomorphism $g$ of $M$ sends circles into 
circles if and only if $g\in G$. More precisely, given a diffeomorphism $g$ 
of $M$, the following assertions are equivalent.

\smallskip

\emph{a)} The diffeomorphism $g$ belongs to $G$.

\smallskip

\emph{b)} For any circle $c$ in $M$, $g\circ c$ is a circle in $M$.

\smallskip

\emph{c)} For any circle $c$ in $M$, $g\circ c$ is a reparametrization of a 
circle in $M$. 
\end{theorem}

\subsection{Diametrical geodesics of symmetric Riemannian metrics on $R$%
-spaces}

\smallskip

It is well-known that any maximal compact subgroup $K$ of $G$ (they are all 
conjugate with each other) acts transitively on $M$ and there is a unique 
(up to homotheties) Riemannian metric $g$ on $M$ such the action of $K$ is 
by isometries of $g$. By abuse of nomenclature we refer to it as \emph{the} $
K$-invariant Riemannian metric on $M$. Moreover, $\left( M,g\right) $ is a 
compact, irreducible Riemannian symmetric space with cubic maximal tori (see 
\cite{Loos}, and also \cite{Jost} for a geometrical proof of the
rectangularity of  maximal tori).

Let $T$ be an $r$-dimensional cubic flat torus of volume $\lambda ^{r}$, 
that is, $T=\mathbb{R}^{r}/\lambda \mathbb{Z}^{r}$ endowed with the 
Riemannian metric such that the canonical projection $p:\mathbb{R}
^{r}\rightarrow T$ is a local isometry. A geodesic $\gamma :\mathbb{R}
\rightarrow T$ is said to be \emph{diagonal }if there is an isometry $F$ of $
T$ such that $\gamma \left( ct\right) =F\left( p\left( t,\dots ,t\right) 
\right) $ for some $c\neq 0$ and all $t$.

Let $M$ be a compact Riemannian manifold with associated distance $d$. The 
diameter of $M$, denoted by diam$~M$, is the maximum of the numbers $d\left(
p,q\right) $ among all pairs of points $p,q$ in $M$. Two points $p,q\in M$ 
are said to be \emph{diametrical} if $d\left( p,q\right) =$ diam~$M$. If $M$
is symmetric, a unit speed geodesic $\gamma :\mathbb{R}\rightarrow M$ is 
said to be \emph{diametrical} if $\gamma \left( t\right) ,\gamma \left( t+ 
\text{diam~}M\right) $ are diametrical points in $M$ for all $t$. A 
nonconstant geodesic in $M$ is said to be diametrical if some (any) unit 
speed reparametrization is so.

As a corollary of the main result in \cite{sakai}, we have that a geodesic
in  a symmetric $R$-space is diametrical if and only if it is a diagonal 
geodesic in a (cubic) maximal torus of $N$, in particular periodic with 
period $2$~diam~$M$ if it is parametrized by arc length.

From now on, we set $\tan \left( \frac{\pi }{2}+k\pi \right) =\infty $ for 
any $k\in \mathbb{Z}$.

\begin{theorem}
\label{Dos}Let $c:\mathbb{R}\cup \left\{ \infty \right\} \rightarrow M$ be a
circle in a self dual symmetric $R$-space $M$ with big transformation group $
G$. Then there exists a maximal compact subgroup $K\ $of $G$ such that the 
curve $\gamma :\mathbb{R}\rightarrow M$ defined by  
\begin{equation}
\gamma \left( t\right) =c\left( \tan \left( \pi t\right) \right)
\label{gammac}
\end{equation}
is a diametrical geodesic of $M$, provided that $M$ is endowed with the 
symmetric $K$-invariant Riemannian metric .

Conversely, given a diametrical geodesic $\gamma $ with period $1$ of a 
Riemannian self dual symmetric $R$-space $M$, then the curve $c:\mathbb{R}
\cup \left\{ \infty \right\} \rightarrow M$ uniquely determined by \emph{(\ref{gammac})}
 is a circle in $M$. 
\end{theorem}

The geodesic and the circle coincide up to a projective reparametrization in
the following sense: Let $f,\varphi $ be the functions 
\begin{equation*}
\mathbb{R}\cup \left\{ \infty \right\} \overset{f}{\longleftarrow }\mathbb{R}
/\mathbb{Z}\overset{\varphi }{\longrightarrow }S^{1}\subset \mathbb{C}
\end{equation*}
defined by $f\left( t+\mathbb{Z}\right) =\tan \left( \pi t\right) $ and $
\varphi \left( t\right) =e^{2\pi ti}$. Then $\varphi \circ f^{-1}$ is the 
projective map $s\mapsto \left( 1-s^{2}+2si\right) /\left( 1+s^{2}\right) $.

\bigskip

I would like to express my gratitude to Cristi\'{a}n U.\ S\'{a}nchez for 
introducing and promoting in C\'{o}rdoba the beautiful and central concept 
of symmetric $R$-space.

\section{Preliminaries}

\subsection{Prevalent vectors}

We fix a self dual symmetric $R$-space $M$, the orbit of height one
parabolic subalgebra of the simple Lie algebra $\mathfrak{g}$ under the
group Inn~$\left( \mathfrak{g}\right) $. The following definition, taken
from \cite{Frost}, is motivated by Lemma 4.1 in \cite{BP} and Proposition %
\ref{equivalencias} below is a restatement of the latter. It provides a Lie
algebraic condition on a tangent vector to $M$ to be the initial velocity of
a circle in $M$.

\medskip

\noindent \textbf{Definition. }Let $M$ be a self dual symmetric $R$-space 
with big transformation group $G$. A vector $y\in \mathfrak{g}$ is said to 
be \emph{prevalent} if  
\begin{equation}
y\in \mathfrak{p}^{\bot }\text{\ \ \ \ \ \ \ \ \ \ and\ \ \ \ \ \ \ \ \ \
Ker }\,\left( \text{ad}_{y}\right) ^{2}=\mathfrak{p}  \label{prevalent}
\end{equation}
for some $\mathfrak{p}\in M$.

\medskip

\noindent \textbf{Remarks. }a) In the second equation the inclusion $\supset
$ is always valid.

\smallskip

b) By Proposition 7.6 in \cite{Frost}, the notion of prevalence makes sense
only for self dual symmetric $R$-spaces. Also, in that Ph.D. thesis,
prevalent vectors are called regular vectors, but we prefer a more specific
term. The importance of the concept was apparent in \cite{BP}, but they did
not give it a name. It was also present in \cite{GK} (the nondegeneracy of $%
P\left( X\right) $), but involving a Cartan involution $\tau $, which is an
alien structure in our setting (cf.\ the proof of Proposition \ref%
{Vprevalente} below).

\medskip

\begin{proposition}
\label{equivalencias}Let $M$ be a self dual symmetric $R$-space with big 
transformation group $G$. Let $\mathfrak{p}$ and $\mathfrak{q}$ be two 
opposite elements of $M$ and let $y\in \mathfrak{q}^{\bot }$. The following 
assertions are equivalent:

\emph{a)} The vector $y$ is prevalent.

\emph{b)} The map 
\begin{equation}
\left. \left( \text{\emph{ad}}_{y}\right) ^{2}\right\vert _{\mathfrak{p}}: 
\mathfrak{p}^{\bot }\rightarrow \mathfrak{q}^{\bot }  \label{iso}
\end{equation}
is a linear isomorphism.

\emph{c)} For any $t\neq 0$, $\exp \left( ty\right) \cdot \mathfrak{p}$ is 
an element of $M$ opposite to $\mathfrak{p}$ and $\mathfrak{q}$.

\emph{d)} The curve $c:\mathbb{R}\cup \left\{ \infty \right\} \rightarrow M$
, $c\left( t\right) =\exp \left( ty\right) \cdot \mathfrak{p}$, $c\left( 
\infty \right) =\mathfrak{q}$, is the circle through $\mathfrak{p}$, $\exp 
\left( y\right) \cdot \mathfrak{p}$ and $\mathfrak{q}$. 
\end{proposition}

\begin{proposition}
Given $y\in \mathfrak{g}$, there is at most one $\mathfrak{p}\in M$ 
satisfying condition \emph{(\ref{prevalent})}. 
\end{proposition}

\noindent \textit{Proof. }Suppose that $\mathfrak{q}\in M$ also satisfies 
conditions (\ref{prevalent}). Since $\mathfrak{q}$ and $\mathfrak{p}$ have 
the same dimension, it suffices to check that $\mathfrak{p}\subset \mathfrak{%
\ q}$. Let $z\in \mathfrak{p}$. Since $y\in \mathfrak{p}^{\bot }$, ad$%
_{y}z\in  \left[ \mathfrak{p}^{\bot },\mathfrak{p}\right] \subset \mathfrak{p%
}^{\bot }$ . Hence, ad$_{y}$~ad$_{y}z\in \left[ \mathfrak{p}^{\bot },%
\mathfrak{p}^{\bot }\right] =\left\{ 0\right\} $ ($\mathfrak{p}^{\bot }$ is
abelian). Thus $ \left( \text{ad}_{y}\right) ^{2}z=0$ and so $z\in $ Ker~$%
\left( \text{ad} _{y}\right) ^{2}=\mathfrak{q}$. \hfill $\square $

\bigskip

\noindent \textbf{Remark. }We wonder whether prevalent vectors in $\mathfrak{%
\ g}$ can be characterized as those elements of $\mathfrak{g}$ belonging to 
exactly one $\mathfrak{p}^{\bot }$ with $\mathfrak{p}\in M$.

\bigskip

The notion of prevalence of a vector in $\mathfrak{g}$ defines the concept
of prevalent tangent vector of $M$:

\smallskip

\noindent \textbf{Definition. }Let $M$ be a self dual symmetric $R$-space. A
tangent vector $X\in T_{\mathfrak{q}}M$ is said to be \emph{prevalent} if $
X=\left. \frac{d}{dt}\right\vert _{0}e^{ty}\cdot \mathfrak{q}$ for some 
prevalent vector $y\in \mathfrak{p}^{\bot }$ with $\mathfrak{p}$ opposite to
$\mathfrak{q}$.

\smallskip

This notion is independent of the choice of $\mathfrak{p}\in \Omega _{ 
\mathfrak{q}}$. Indeed, suppose that $\mathfrak{p}^{\prime }\in M$ is 
opposite to $\mathfrak{q}$. Then there exists a unique $z\in \mathfrak{q}
^{\bot }$ such that $\mathfrak{p}^{\prime }=e^{z}\cdot \mathfrak{p}$. We 
have that  
\begin{equation*}
e^{z}y=y+\left[ z,y\right] +\frac{1}{2}\left[ z,\left[ z,y\right] \right]
\end{equation*}
(see the proof of Lemma 4.1 in \cite{BP}), where the second and third terms 
of the right hand side are in $\left[ \mathfrak{q}^{\bot },\mathfrak{p}
^{\bot }\right] \subset \mathfrak{q}\cap \mathfrak{p}$ and $\left[ \mathfrak{%
\ q}^{\bot },\mathfrak{q}\right] $, respectively, and hence both are
contained  in $\mathfrak{q}$. Let $\pi :G\rightarrow M$ be the projection
given by $\pi  \left( g\right) =g\cdot \mathfrak{q}$, which satisfies that
Ker~$\left( d\pi  _{e}\right) =\mathfrak{q}$. Hence,  
\begin{equation*}
X=d\pi _{e}\left( y\right) =d\pi _{e}\left( y+\left[ z,y\right] +\frac{1}{2} %
\left[ z,\left[ z,y\right] \right] \right) =\left. \frac{d}{dt}\right\vert
_{0}\exp \left( te^{z}y\right) \cdot \mathfrak{q}\text{,}
\end{equation*}
where $e^{z}y\in \mathfrak{p}^{\prime }$ is prevalent, since it is conjugate
to $y$.

\bigskip

For further reference we recall that given opposite elements $\mathfrak{p}$ 
and $\mathfrak{q}$ in a symmetric $R$-space $M$, there exists $z\in  
\mathfrak{g}$, called the characteristic element of the pair $\mathfrak{p}, 
\mathfrak{q}$, satisfying that $\mathfrak{p}^{\bot },\mathfrak{p}\cap\mathfrak{q}$ 
and $\mathfrak{q}^{\bot }$ are the eigenspaces of the operator ad$
_{z}$ with eigenvalues $1,0$ and $-1$, respectively.

\subsection{Generalized conformal structures}

We collect some facts from \cite{GK} and \cite{libro} Part II on 
characterizations of the big transformation group of a symmetric $R$-space.

Let $V$ be a finite dimensional vector space. An algebraic subset $C$ of $V$
is called a prehomogeneous generalized cone if it is stable by homotheties
and the subgroup of $GL\left( V\right) $ preserving $C$ has an open orbit.
Let $M$ be a smooth manifold with the same dimension as $V$. A \emph{locally
flat generalized conformal structure} on $M$ with typical cone $C$ assigns
to each point $q\in M$ a subset $C_{q}$ of $T_{q}M$ in such a way that for
any $p\in M$ there is an open neighborhood $U$ of $p$ and a smooth map $%
\varphi :U\times V\rightarrow TM$ such that for each $q\in U$, $\varphi
\left( q,\cdot \right) $ is a linear isomorphism from $V$ to $T_{\mathfrak{q}%
}M$ sending $C$ to $C_{q}$.

Given a symmetric $R$-space $M$, Gindikin and Kaneyuki define in Lemmas 3.1 
and 3.2 of \cite{GK} a locally flat generalized conformal structure $ 
\mathcal{K}$ on $M$, where $C=\partial V_{r}$ for a certain subset $V_{r}$ 
of $\mathfrak{q}^{\bot }$ for some $\mathfrak{q}\in M$ ($\partial $ denotes 
the frontier), and prove the following theorem (Theorem 3.3).

\begin{theorem}
\label{mainGK}\emph{\cite{GK}} The group of automorphisms of $\mathcal{K}$ 
is exactly $G$, provided that the symmetric $R$-space $M$ has rank larger 
than one. 
\end{theorem}

In the case of a symmetric $R$-space $M$ of type $C$ (or equivalently, if $M$
is self dual, see the proof of Proposition \ref{Vprevalente} below), 
essentially the same result has been obtained by Bertram in \cite{Bertram}.

\begin{proposition}
\label{Vprevalente}Let $M$ be a self dual symmetric $R$-space. Then, for all
$\mathfrak{p}\in M$,  
\begin{equation*}
\mathcal{K}_{\mathfrak{p}}=\partial \left\{ X\in T_{\mathfrak{p}}M\mid X 
\text{ is prevalent}\right\} \text{.}
\end{equation*}
\end{proposition}

\noindent \textit{Proof. }First, we recall the definition of the generalized
conformal structure $\mathcal{K}$ on $M$:\ Given $\mathfrak{p}\in M$, take 
any $\mathfrak{q}\in M$ opposite to $\mathfrak{p}$ and define $\mathcal{K}_{ 
\mathfrak{p}}=\pi _{\mathfrak{q}}^{\mathfrak{p}}\left( \partial V_{r}\right)
,$ where $V_{r}$ turns out to be equal to $\left\{ x\in \mathfrak{q}^{\bot
}\mid x\text{ is prevalent}\right\} $, provided that $M$ is self dual. More 
precisely, let $z$ be the characteristic element of the pair $\mathfrak{p}, 
\mathfrak{q}$ and define $P:\mathfrak{q}^{\bot }\rightarrow $ End~$\left(  
\mathfrak{q}^{\bot }\right) $ as in 2.7 of \cite{GK} by  
\begin{equation*}
P\left( x\right) \left( y\right) =\tfrac{1}{2}\left[ \left[ \tau y,x\right]
,x\right] =\tfrac{1}{2}\left( \text{ad}_{x}\right) ^{2}\left( \tau y\right) 
\text{,}
\end{equation*}
where $\tau $ is a Cartan involution of $\mathfrak{g}$ satisfying $\tau 
\left( z\right) =-z$.

Now, comparing the table of self dual symmetric $R$-spaces at the end of  
\cite{BP} with Table 2 (page 601) in \cite{KaneSylves} (or Table 5 on page 
135 of \cite{libro} Part II), we see that a symmetric $R$-space is self dual
if and only if certain root system $\Delta \left( \mathfrak{g},\mathfrak{c}
\right) $ is of type $C$ (there should be a better reason for this to be 
true). In this case, one has by (2.24) in \cite{GK} that  
\begin{equation}
V_{r}=\left\{ x\in \mathfrak{q}^{\bot }\mid \det P\left( x\right) \neq
0\right\} .  \label{cone}
\end{equation}
By the definition of prevalent tangent vectors in $T_{\mathfrak{p}}M$, it 
suffices to check that $V_{r}$ consists exactly of the prevalent vectors in $
\mathfrak{q}^{\bot }$. Indeed, suppose that $x\in V_{r}$, that is, $\left(  
\text{ad}_{x}\right) ^{2}\circ \tau :\mathfrak{q}^{\bot }\rightarrow  
\mathfrak{q}^{\bot }$ is a linear isomorphism. But since $\left. \tau
\right\vert  _{\mathfrak{q}^{\bot }}:\mathfrak{q}^{\bot }\rightarrow 
\mathfrak{p}^{\bot }$  is an isomorphism, then $\left( \text{ad}_{x}\right)
^{2}:\mathfrak{p}^{\bot }\rightarrow \mathfrak{q}^{\bot }$ also is an
isomorphism. Thus, $x$ is  prevalent. The other inclusion holds by similar
arguments. \hfill $\square $

\section{Proofs of the theorems}

\noindent \textit{Proof of Theorem \ref{Uno}. }First we check that a) $
\Rightarrow $ b). Let $\mathfrak{p},\mathfrak{p}_{1}$ and $\mathfrak{q}$ be 
three points in $M$ which are pairwise opposite and let $c:\mathbb{R}\cup 
\left\{ \infty \right\} \rightarrow M$ be the circle through $\mathfrak{p}, 
\mathfrak{p}_{1},\mathfrak{q}$, that is, $c\left( t\right) =e^{ty}\cdot  
\mathfrak{p}$ and $c\left( \infty \right) =\mathfrak{q}$ for a unique $y\in  
\mathfrak{q}^{\bot }$. Let $g\in G$. Since $g$ is an automorphism of $ 
\mathfrak{g}$, then $g\cdot \mathfrak{p},g\cdot \mathfrak{p}_{1},g\cdot  
\mathfrak{q}$ are pairwise opposite in $M$ and $g$ also preserves the 
Killing form; in particular, $g\cdot y\in $ $g\cdot \mathfrak{q}^{\bot
}=\left( g\cdot \mathfrak{q}\right) ^{\bot }$. Hence,  
\begin{equation*}
g\left( c\left( t\right) \right) =ge^{ty}\cdot \mathfrak{p}=e^{tgy}\cdot g 
\mathfrak{p}\text{, \ \ \ \ \ \ \ \ \ \ \ }g\left( \infty \right) =g\cdot 
\mathfrak{q}
\end{equation*}
is the circle through $g\cdot \mathfrak{p},g\cdot \mathfrak{p}_{1},g\cdot  
\mathfrak{q}$.

Since the implication b) $\Rightarrow $ c) is obvious, it remains only to
prove c) $\Rightarrow $ a). Suppose that $g$ is a diffeomorphism of $M$
satisfying c). For the sphere $S^{n}$ we know from \cite{Jeffers} that the
assertion is true (see also Theorem 3 in \cite{CPV}, a much more general
result). The projective spaces $P^{n}\left( \mathbb{F}\right) $, where $%
\mathbb{F}=\mathbb{C},\mathbb{H}$ or $\mathbb{O}$ ($n=2$ in the last case),
are not self dual (see the table at the end of \cite{BP}). Thus, we may
suppose that rank~$\left( M\right) >1$. By Theorem \ref{mainGK}, it suffices
to prove that $g$ is an automorphism of the generalized conformal structure $%
\mathcal{K}$.

Given $\mathfrak{p}\in M$, we must show that $dg_{\mathfrak{p}}\left( 
\mathcal{K}_{\mathfrak{p}}\right) =\mathcal{K}_{g\left( \mathfrak{p}\right) }
$. Since $dg_{\mathfrak{p}}$ is a linear isomorphism, by Proposition \ref%
{Vprevalente} it is equivalent to show that $dg_{\mathfrak{p}}\left(
X\right) $ is prevalent for any prevalent tangent vector $X\in T_{\mathfrak{p%
}}M$. Suppose $X=\pi _{\mathfrak{q}}^{\mathfrak{p}}\left( x\right) $ for
some prevalent vector $x\in \mathfrak{q}^{\bot }$. By Proposition \ref%
{equivalencias}, the curve $c:\mathbb{R}\cup \left\{ \infty \right\}
\rightarrow M$, $c\left( t\right) =\exp \left( tx\right) \cdot \mathfrak{p}$
, $c\left( \infty \right) =\mathfrak{q}$ is a circle in $M$. Now, by the
hypothesis, $g\circ c=\bar{c}\circ \phi $ for some circle $\bar{c}$ in $M$
and some diffeomorphism $\phi $ of $\mathbb{R}\cup \left\{ \infty \right\} $
. Suppose that $\bar{c}$ is the circle through $\overline{\mathfrak{p}}$, $%
\overline{\mathfrak{p}}_{1}$ and $\mathfrak{\bar{q}}$, that is, $\bar{c}%
\left( t\right) =\exp \left( t\overline{x}\right) \cdot \mathfrak{\bar{p}}$
for some $\overline{x}\in \mathfrak{\bar{q}}^{\bot }$ with $\mathfrak{\bar{q}%
}$ opposite to $\overline{\mathfrak{p}}$ and $\overline{\mathfrak{p}}_{1}$.
In particular, $\overline{x}$ is prevalent by Proposition \ref{equivalencias}.
 We call $t_{o}=\phi \left( 0\right) $ and assuming that $t_{o}\neq \infty $,
 we compute 
\begin{equation*}
dg_{\mathfrak{p}}\left( X\right) =\left( g\circ c\right) ^{\prime }\left(
0\right) =\left( \bar{c}\circ \phi \right) ^{\prime }\left( 0\right) =\bar{c}%
^{\prime }\left( t_{o}\right) \phi ^{\prime }\left( 0\right) \text{,}
\end{equation*}%
with $\phi ^{\prime }\left( 0\right) \neq 0$. Now, calling $h=\exp \left(
t_{o}\overline{x}\right) $, we have that 
\begin{equation*}
\bar{c}^{\prime }\left( t_{o}\right) =\left. \frac{d}{ds}\right\vert _{s=0}%
\bar{c}\left( t_{o}+s\right) =\left. \frac{d}{ds}\right\vert _{s=0}\exp
\left( t_{o}\overline{x}\right) \exp \left( s\overline{x}\right) \cdot 
\overline{\mathfrak{p}}=\left. \frac{d}{ds}\right\vert _{s=0}\exp \left( sh%
\overline{x}\right) \cdot \overline{\mathfrak{p}}=\pi _{h\mathfrak{\bar{q}}%
}^{h\overline{\mathfrak{p}}}\left( h\overline{x}\right) \text{.}
\end{equation*}%
Since $h$ is an automorphism and $\overline{x}$ is prevalent, we have that $%
dg_{\mathfrak{p}}\left( X\right) $ is prevalent. Therefore, $g$ is an
automorphism of $\mathcal{K}$, as desired. For the case $t_{o}=\infty $ use
a reparametrization as in the last assertion of Subsection \ref{SScircles}.
\hfill $\square $

\bigskip

\noindent \textit{Proof of Theorem \ref{Dos}. }The general setting of the 
proof is from \cite{BasicTakeuchi} and it uses results from \cite{GK,
KaneSylves}. Only we should keep in mind that in contrast to those articles, 
no maximal compact subgroup of $G$ (or equivalently, no symmetric Riemannian
metric on $M$) is fixed. The notation is a mixture of that of the different 
sources. Suppose that the circle $c$ is given by $c\left( 0\right) = 
\mathfrak{p}$, $c\left( \infty \right) =\mathfrak{q}$ and $c\left( t\right) 
=\exp \left( tx\right) \cdot \mathfrak{p}$ for all $t$, with $x\in \mathfrak{%
\ q}^{\bot }$; in particular, $x$ is prevalent.

We consider the decomposition $\mathfrak{g}=\mathfrak{p}^{\bot }\oplus  
\mathfrak{g}_{0}\oplus \mathfrak{q}^{\bot }$, with $\mathfrak{g}_{0}= 
\mathfrak{p}\cap \mathfrak{q}$, which turns out to be the Lie algebra of $
G_{0}\left( \mathfrak{p},\mathfrak{q}\right) $. Let $z\in \mathfrak{g}_{0}$ 
be the associated characteristic element. By Theorem I.2.3 in  \cite{libro}
Part II there exists a Cartan involution $\tau $ on $\mathfrak{g}$ such 
that $\tau \left( \mathfrak{p}^{\bot }\right) =\mathfrak{q}^{\bot }$ and $
\tau \left( \mathfrak{g}_{0}\right) =\mathfrak{g}_{0}$; in particular, $\tau
\left( z\right) =-z$ (see (1.4) in \cite{BasicTakeuchi}). Let $\mathfrak{g}= 
\mathfrak{k}\oplus\mathfrak{n}$ be the decomposition associated with $\tau $, 
that is, $\mathfrak{k}$ and $\mathfrak{n}$ are the eigenspaces of $\tau $ 
with eigenvalues 1 and $-1$. Define 
\begin{equation}
\overline{K}=\left\{ g\in G\mid g\tau =\tau g\right\} \text{,}  \label{Kbar}
\end{equation}
which is a maximal compact subgroup of $G$, with Lie algebra $\mathfrak{k}$,
acting transitively on $M$, with isotropy subgroup $\overline{K}_{0}=_{\text{
def}}\overline{K}\cap G_{0}\left( \mathfrak{p},\mathfrak{q}\right) $ at $ 
\mathfrak{p}$. Thus, one can identify $M=\overline{K}/\overline{K}_{0}$.

Let $\mathfrak{a}$ be a maximal abelian subspace of $\mathfrak{n}$ 
containing $z$. Let $\Sigma =\Delta \left( \mathfrak{g},\mathfrak{a}\right) $
be the irreducible root system of $\mathfrak{g}$ relative to $\mathfrak{a}$ 
and $\Sigma _{i}=\left\{ \gamma \in \Sigma \mid B\left( \gamma ,z\right) 
=i\right\} $ for $i=0,\pm 1$, where $B$ denotes the Killing form of $ 
\mathfrak{g}$.

Choose a linear order in $\Sigma $ such that the set of positive roots $
\Sigma ^{+}$ satisfies $\Sigma _{0}\subset \Sigma ^{+}\subset \Sigma 
_{0}\cup \Sigma _{1}$ and a maximal system of strongly orthogonal roots $
\left\{ \beta _{1},\dots ,\beta _{r}\right\} $ in $\Sigma _{1}$ such that $
\beta _{1}$ is the highest root in $\Sigma $. Let $\mathfrak{g}^{\alpha }$ 
be the root space for the root $\alpha \in \Sigma $. For each $i=1,\dots ,r$
choose a vector $X_{i}\in \mathfrak{g}^{\beta _{i}}\subset \mathfrak{p}
^{\bot }$ with $\left\vert X_{i}\right\vert ^{2}=2/B\left( \beta _{i},\beta 
_{i}\right) $ and set $X_{-i}=\tau \left( X_{i}\right) \in \mathfrak{g}
^{-\beta }\subset \mathfrak{q}^{\bot }$. Then one has $\left[ X_{\beta
_{i}},X_{\beta _{i}}\right] =-2\beta _{i}/B\left( \beta _{i},\beta 
_{i}\right) $.

Now, $\mathfrak{k}$ decomposes as $\mathfrak{k}=\mathfrak{k}_{0}\oplus\mathfrak{m}
_{\mathfrak{k}}$, where $\mathfrak{k}_{0}=\mathfrak{g}_{0}\cap \mathfrak{k}$
and $\mathfrak{m}_{\mathfrak{k}}=\left( \mathfrak{p}^{\bot }\oplus\mathfrak{q}
^{\bot }\right) \cap \mathfrak{k}$. Via the differential of the canonical 
projection $\overline{K}\rightarrow \overline{K}/\overline{K}_{0}$ one 
identifies $T_{o}M=\mathfrak{m}_{\mathfrak{k}}$. Let $A_{i}=\pi \left( 
X_{i}+X_{-i}\right) \in \mathfrak{m}_{\mathfrak{k}}$. By (1.21) in \cite%
{BasicTakeuchi}, $A_{1},\dots ,A_{r}$ is a basis of a maximal abelian 
subspace of $\mathfrak{m}_{\mathfrak{k}}$ and moreover its $\mathbb{Z}$-span
is the lattice of a cubic maximal torus $T$ in $M$. Let $A=
\sum_{i=1}^{r}A_{i}$. Then $\bar{\gamma}\left( t\right) =_{\text{def}}\exp 
\left( tA\right) \overline{K}_{0}$ is a diagonal geodesic in $T$, and hence 
diametrical in $M$.

From the proof of Lemma 2.3 in \cite{BasicTakeuchi} we have that 
\begin{equation*}
\exp \left( \sum_{i=1}^{r}x_{i}X_{-i}\right) \cdot \mathfrak{p}=\exp \left( 
\frac{1}{\pi }\sum_{i=1}^{r}\arctan \left( x_{i}\right) A_{i}\right) \cdot 
\mathfrak{p}
\end{equation*}
for any $x_{1},\dots ,x_{n}\in \mathbb{R}$. In particular, taking $
x_{i}=\tan \left( \pi t\right) $ for all $i$ and calling $\overline{X}
=\sum_{i=1}^{r}X_{-i}\in \mathfrak{q}^{\bot }$, we have  
\begin{equation*}
\bar{\gamma}\left( t\right) =\bar{c}\left( \tan \left( \pi t\right) \right) 
\text{,}
\end{equation*}
where $\bar{c}\left( s\right) =\exp \left( s\overline{X}\right) \cdot  
\mathfrak{p}$. Now, $\overline{X}=\pm O_{r,0}$ in the notation of (1.17) in  
\cite{GK} and by Theorem 2.9 there, it belongs to $V_{r}$. 
Hence, by the proof of Proposition \ref{Vprevalente}, $\overline{X}$ is 
prevalent and so $\bar{c}$ is a circle in $M$ by Proposition \ref%
{equivalencias}.

On the other hand, by \cite{KaneSylves}, the $G_{0}\left( \mathfrak{p}, 
\mathfrak{q}\right) $-orbit of any element in $V_{r}$ is the whole $V_{r}$ 
(actually, Kaneyuki studies the finer $G_{0}^{0}\left( \mathfrak{p}, 
\mathfrak{q}\right) $-orbit structure, where the superscript $0$ refers to 
the identity component, but we do not need that). Hence, there is $h\in 
G_{0}\left( \mathfrak{p},\mathfrak{q}\right) $ such that $h\overline{X}=X$. 
Let us show that $K=_{\text{def}}h^{-1}\overline{K}h$ is a maximal compact 
subgroup of $G$ satisfying the required property. Indeed, since $h\mathfrak{%
p }=\mathfrak{p}$, we have that 
\begin{equation*}
h\bar{c}\left( s\right) =h\exp \left( s\overline{X}\right) \cdot \mathfrak{p}
=h\exp \left( s\overline{X}\right) h^{-1}h\cdot \mathfrak{p}=\exp \left( sh 
\overline{X}\right) \cdot \mathfrak{p}=\exp \left( sX\right) \cdot \mathfrak{%
\ p}=c\left( s\right)
\end{equation*}
for all $s$. Hence  
\begin{equation*}
\gamma \left( t\right) =c\left( \tan \left( \pi t\right) \right) =h\bar{c}
\left( \tan \left( \pi t\right) \right) =h\bar{\gamma}\left( t\right)
\end{equation*}
holds for all $t$. Call $\bar{g}$ the $\overline{K}$-invariant Riemannian 
metric on $M$. Hence $g=h^{\ast }\bar{g}$ is a metric in $M$ which is $K$
-invariant. Since $h:\left( M,\bar{g}\right) \rightarrow \left( M,g\right) $
is an isometry and $\bar{\gamma}$ is a diametrical geodesic of $\left( M, 
\bar{g}\right) $, then the same is true for $\gamma =h\circ \bar{\gamma}$ in
$\left( M,g\right) $, as desired.

Now we deal with the converse. The fact that $M$ has a fixed Riemannian 
structure of that type is equivalent to fixing a subgroup $\overline{K}$ as 
in (\ref{Kbar}). One follows the arguments above until the conclusion that $ 
\bar{c}$ is a circle (since the maximal torus is cubic, one can consider any
diagonal).\hfill $\square $

\section{Examples\label{Examples}}

The examples below complement those in Subsection 2.2 of \cite{BP} or 
approach them from a more concrete (but also quite general) perspective; for
instance, the seven families of split isotropic Grassmannians and their 
circles are treated in a uniform way. We include the circles of the 
Grassmannian of oriented planes (aka complex quadric), which is a self dual
symmetric $R$-space not described in detail in \cite{BP} and is very 
different from the split standard or isotropic Grassmannians (see Subsection
\ref{quadric} below).

In each of the following examples the smooth manifold $M$ is acted on 
transitively by a semisimple Lie group $H$. Each point $p\in M$ is 
identified with its \emph{infinitesimal stabilizer} $\mathfrak{p}$, that is,
the Lie algebra of the isotropy subgroup of $H$ at $\mathfrak{p}$. If $H$ is
simple, then $M$ is a symmetric $R$-space, $H$ coincides with the big 
transformation group $G$ of $M$ up to coverings and connected components and
$\mathfrak{p}$ is a height one parabolic subalgebra of $\mathfrak{g}=$ Lie~$
\left( G\right) =$ Lie~$\left( H\right) $.

\subsection{The sphere}

The group $H=O_{o}\left( 1,n+1\right) $ acts on the sphere $S^{n}$ by 
orientation preserving conformal diffeomorphisms. Let $\mathfrak{q}$ be the 
infinitesimal stabilizer of a point $q\in S^{n}$ and let $
F_{q}:S^{n}\rightarrow q^{\bot }\cup \left\{ \infty \right\} $ be the 
standard stereographic projection with $F_{q}\left( q\right) =\infty $. 
Then, after conjugation by $F_{q}$, $\exp \left( \mathfrak{q}^{\bot }\right)
$ consists of the diffeomorphisms of $q^{\bot }\cup \left\{ \infty \right\} $
fixing $\infty $ and acting on $q^{\bot }$ by translations. Two points in $
S^{n}$ are opposite if and only if they are different. The trajectories of 
the circles are the usual circles: intersections of $S^{n}$ with planes at 
distance smaller than one from the origin. The circle $c$ determined by 
three distinct points $p,p_{1}$ and $q$ is given by $F_{q}\left( c\left( 
t\right) \right) =\left( 1-t\right) F_{q}\left( p\right) +tF_{q}\left( 
p_{1}\right) $ and $c\left( \infty \right) =q$.

\subsection{The split standard Grassmannians\label{ssg}}

Let $\mathbb{F}=\mathbb{R},\mathbb{C}$ or $\mathbb{H}$ (where the skew field
$\mathbb{H}$ of the quaternions acts on $\mathbb{H}^{m}$ on the right). The 
group $H=SL\left( m,\mathbb{F}\right) $ acts in a canonical way on the 
Grassmannian $M$ of $k$-dimensional subspaces of $\mathbb{F}^{m}$, which is 
a symmetric $R$-space. A subspace $P$ is identified with the subalgebra $ 
\mathfrak{p}$ of $sl\left( m,\mathbb{F}\right) $ consisting of all the $ 
\mathbb{F}$-linear transformations of $\mathbb{F}^{m}$ preserving $P$ and $ 
\mathfrak{p}^{\bot }$ consists of the $\mathbb{F}$-linear transformations 
vanishing on $P$ with image contained in $P$. Two subspaces $P,Q\in M$ are 
opposite if and only if $P\cap Q=\left\{ 0\right\} $. Now, $M$ is self dual 
if and only if $m=2n$. In this case, which is the only we consider 
henceforth, circles in $M$ have the form 
\begin{equation*}
c\left( t\right) =\left\{ x+tT\left( x\right) \mid x\in P\right\} \text{, \
\ \ \ \ \ \ \ \ \ \ }c\left( \infty \right) =Q\text{,}
\end{equation*}
for some complementary subspaces $P,Q$ of the same dimension in $M$ and some
$\mathbb{F}$-linear isomorphism $T:P\rightarrow Q$. Let $\left\{ u_{1},\dots
,u_{n}\right\} $ be an $\mathbb{F}$-basis of $P$ and let $g$ be the $\mathbb{%
\ F}$-Hermitian symmetric (positive definite) inner product on $\mathbb{F}
^{2n} $ such that $\left\{ u_{1},\dots ,u_{n},Tu_{1},\dots ,Tu_{n}\right\} $
is orthonormal. Then $K=_{\text{def}}~$Aut$_{0}$~$\left( g\right) $ is a 
maximal compact subgroup of $H$ and $\gamma :\mathbb{R}\rightarrow M$ 
defined by  
\begin{equation*}
\gamma \left( s\right) =\text{span}~\left\{ \cos \left( \pi s\right)
u_{i}+\sin \left( \pi s\right) T\left( u_{i}\right) \mid i=1,\dots ,n\right\}
\end{equation*}
is a geodesic in $M$ provided that it carries the (unique up to homotheties)
$K$-invariant Riemannian metric. Clearly, (\ref{gammac}) is satisfied.

\subsection{The split isotropic Grassmannians associated with the seven
split inner products\label{sig}}

We recall some basic facts about bilinear and sesquilinear forms from 
Thorbergsson \cite{Thor} (see also \cite{Harvey, Neretin}). Let $V$ denote a
right vector space over $\mathbb{F}$ where $\mathbb{F}$ is $\mathbb{R}$, $ 
\mathbb{C}$ or $\mathbb{H}$. Let $\sigma :\mathbb{F}\rightarrow \mathbb{F}$ 
be either the identity or conjugation in $\mathbb{F}$, denoted as usual conj$
\ \left( \alpha \right) =\bar{\alpha}$ (in particular, $\sigma $ is the 
identity if $\mathbb{F}=\mathbb{R}$). We recall that $\overline{zw}= 
\overline{w}~\overline{z}$ holds in $\mathbb{H}$.

Let $f:V\times V\rightarrow \mathbb{F}$ be a map that is additive in both 
arguments and nondegenerate, and let $\varepsilon $ be either $1$ or $-1$. 
The map $f$ is said to be an $\varepsilon $-symmetric $\sigma $-sesquilinear
form on $V$, or briefly a $\left( \sigma ,\varepsilon \right) $-form on $V$,
if  
\begin{equation*}
f(x\alpha ,y\beta )=\sigma \left( \alpha \right) f(x,y)\beta \text{\ \ \ \ \
\ \ \ \ \ \ \ \ \ and \ \ \ \ \ \ \ \ \ \ \ \ \ \ }f(x,y)=\varepsilon \sigma
(f(y,x))
\end{equation*}
for all $x$ and $y$ in $V$ and all $\alpha $ and $\beta $ in $\mathbb{F}$. 
Since $\mathbb{H}$ is not commutative, $\sigma $ must be the conjugation if $
\mathbb{F}=\mathbb{H}$. In the case $\mathbb{F}=\mathbb{C}$, $\left( \sigma
,1\right)  $-forms differ inessentially from $\left( \sigma ,-1\right) $%
-forms: If  $f$ is a $\left( \sigma ,-1\right) $-form, then $ if$ is a $%
\left( \sigma ,1\right) $-form.

Let $f$ be a $\left( \sigma ,\varepsilon \right) $-form on $V$. A subspace $
W $ in $V$ is said to be totally isotropic if $f(x,y)=0$ for all $x$ and $y$
in $W$. We will consider only \emph{split} $\left( \sigma ,\varepsilon 
\right) $-forms, that is, the dimension of any maximal totally isotropic 
subspace is half the dimension of $V$; in particular, the dimension of $V$ 
is even. There are seven non-isomorphic families of them. In the table below
we recall the standard models $f:\mathbb{F}^{2n}\times \mathbb{F}
^{2n}\rightarrow \mathbb{F}$. In all cases, $\mathbb{F}^{n}\times \left\{ 
0\right\} $ and $\left\{ 0\right\} \times \mathbb{F}^{n}$ are examples of
maximal  totally isotropic subspaces. 
\begin{equation*}
\begin{tabular}{|c|c|c|c|c|c|}
\hline
$\sigma $ & $\varepsilon $ & $\mathbb{F}$ & Type & $f\left( \left(
x,y\right) ,\left( x^{\prime },y^{\prime }\right) \right) $ & $\mathcal{U}
\left( f\right) $ \\ \hline
id & $1$ & $\mathbb{C}$ & $\mathbb{C}$-symmetric & $x^{\bot }y^{\prime
}+y^{\bot }x^{\prime }$ & $O_{2n}\left( \mathbb{C}\right) $ \\ \hline
id & $-1$ & $\mathbb{R},\mathbb{C}$ & $\mathbb{F}$-symplectic & $x^{\bot
}y^{\prime }-y^{\bot }x^{\prime }$ & $Sp_{2n}\left( \mathbb{F}\right) $ \\ 
\hline
conj & $1$ & $\mathbb{R},\mathbb{C},\mathbb{H}$ & split $\mathbb{F}$
-Hermitian & $\bar{x}^{\bot }y^{\prime }+\bar{y}^{\bot }x^{\prime }$ & $%
O_{n,n}$, $U_{n,n}$, $Sp_{n,n}$ \\ \hline
conj & $-1$ & $\mathbb{H}$ & $\mathbb{H}$-skew Hermitian & $\bar{x}^{\bot
}y^{\prime }-\bar{y}^{\bot }x^{\prime }$ & $SO_{2n}^{\ast }$ \\ \hline
\end{tabular}
\label{table}
\end{equation*}

If $f$ is a split $\left( \sigma ,\varepsilon \right) $-form over $\mathbb{F}
$, we call $M=$ Gr$_{0}$~$\left( f\right) $ the set of all maximal isotropic
subspaces of $V$ for $f$ and put $H=\mathcal{U}\left( f\right) $, the Lie 
group of automorphisms of $f$, which acts transitively on Gr$_{0}$~$\left( 
f\right) $.

All the split isotropic Grassmannians Gr$_{0}$~$\left( f\right) $ are self 
dual, except in the case of a split real Hermitian (o equivalent, split real
symmetric) form on $\mathbb{R}^{2n}$ with $n$ odd (see the table at the end 
of \cite{BP}). A maximal isotropic subspace $P$ is identified with the Lie 
algebra $\mathfrak{p}$ of the parabolic subgroup of \/$\mathcal{U}\left( 
f\right) $ fixing $P$ as a subspace of $V$. Two parabolic subalgebras $ 
\mathfrak{p},\mathfrak{q}$ in $M$ are opposite if and only if the 
corresponding maximal isotropic subspaces $P,Q$ are complementary.

One says that an additive map $S:V\rightarrow V$ is id-linear if it is $ 
\mathbb{F}$-linear and conj-linear (also called $\mathbb{F}$-antilinear) if $
S\left( u\alpha \right) =S\left( u\right) \bar{\alpha}$ for any $u\in V$ and
$\alpha \in \mathbb{F}$. Antilinear maps on linear spaces over $\mathbb{H}$ 
do not exist.

\bigskip

\noindent \textbf{Definition. } Let $f:V\times V\rightarrow \mathbb{F}$ be a
split $\left( \sigma  ,\varepsilon \right) $-form over $\mathbb{F}$. A $\bar{%
\sigma}$-linear map $ J:V\rightarrow V$ is said to be \emph{compactly
adapted to} $f$ if  
\begin{equation}
J^{2}=\varepsilon ~\text{id}\ \ \ \ \ \ \ \ \ \ \ \ \ \ \ and\ \ \ \ \ \ \ \
\ \ \ \ f\left( Jx,Jy\right) =\bar{\sigma}f\left( x,y\right)  \label{ca}
\end{equation}
for all $x,y\in V$ and the symmetric Hermitian form $D:V\times V\rightarrow  
\mathbb{F}$ given by $D\left( x,y\right) =f\left( Jx,y\right) $ is (positive
or negative) definite.

\medskip

In fact, $D$ is a symmetric Hermitian form: We have that 
\begin{equation*}
D\left( x\alpha ,y\beta \right) =f\left( J\left( x\alpha \right) ,y\beta
\right) =f\left( J\left( x\right) \sigma \left( \bar{\alpha}\right) ,y\beta
\right) =\sigma ^{2}\left( \bar{\alpha}\right) f\left( Jx,y\right) \beta = 
\bar{\alpha}D\left( x,y\right) \beta
\end{equation*}
and similarly, $D\left( y,x\right) =\overline{D\left( x,y\right) }$ for all $
x,y\in V$.

\medskip

The following result of Ju.\ Neretin will be helpful to find an appropriate 
symmetric Riemannian metric on Gr$_{0}$~$\left( f\right) $ such that a given
circle is a diametrical geodesic, up to projective reparametrizations.

\begin{lemma}
\label{Neretin}\emph{\cite{Neretin}} Let $f:V\times V\rightarrow \mathbb{F}$
be a split $\left( \sigma ,\varepsilon \right) $-form over $\mathbb{F}$. If $
J:V\rightarrow V$ is compactly adapted to $f$, then $\mathcal{U}^{J}\left( 
f\right) =_{\text{\emph{def}}}\left\{ A\in \mathcal{U}\left( f\right) \mid 
A\circ J=J\circ A\right\} $ coincides with $\mathcal{U}\left( f\right) \cap  
\mathcal{U}\left( D\right) $ and is a maximal compact subgroup of $\mathcal{%
U }\left( f\right) $. 
\end{lemma}

\begin{lemma}
\label{remark}\emph{a)} If $P$ and $Q$ are complementary maximal isotropic 
subspaces of $V$ and $\mathcal{B}=\left\{ u_{1},\dots ,u_{n}\right\} $ is a 
basis of $P$ as an $\mathbb{F}$-space, then there exists a basis $\mathcal{B}
^{\prime }=\left\{ v_{1},\dots ,v_{n}\right\} $ of $Q$, called a dual basis 
of $\mathcal{B}$ with respect to $f$, such that $f\left( u_{i},v_{j}\right) 
=\delta _{ij}$ for all $i,j=1,\dots ,n$.

\smallskip

\emph{b)} Given bases $\mathcal{B}$ and $\mathcal{B}^{\prime }$ of $V$ as 
above, the $\bar{\sigma}$-linear map $J$ such that $Ju_{i}=v_{i}$ and $
Jv_{i}=\varepsilon u_{i}$ is compactly adapted to $f$. 
\end{lemma}

\noindent \textit{Proof. }The first assertion is well-known. Regarding the 
second statement, straightforward computations show that (\ref{ca})  holds
for all $x,y\in V$. We check that the associated form $D$ is $ \varepsilon $%
-definite. Indeed, 
\begin{equation*}
D\left( u_{i},u_{j}\right) =f\left( Ju_{i},u_{j}\right) =f\left(
v_{i},u_{j}\right) =\varepsilon \sigma (f(u_{j},v_{i}))=\varepsilon \delta
_{ij}\text{,}
\end{equation*}
and in the same manner, $D\left( v_{i},v_{j}\right) =\varepsilon \delta 
_{ij} $. \hfill $\square $

\begin{theorem}
Let $f$ be a split $\left( \sigma ,\varepsilon \right) $-form on a vector
space $V$ over $\mathbb{F}$. Let $P$ and $Q$ be two complementary maximal
isotropic subspaces of $V$ and let $T:P\rightarrow Q$ be an $\mathbb{F}$
-linear isomorphism satisfying 
\begin{equation}
f\left( Tx,y\right) +f\left( x,Ty\right) =0  \label{Tskew}
\end{equation}%
for all $x,y\in P$. Then $c:\mathbb{R}\cup \left\{ \infty \right\}
\rightarrow $ Gr$_{0}$~$\left( f\right) $ defined by 
\begin{equation}
c\left( t\right) =\left\{ x+tT\left( x\right) \mid x\in P\right\} \text{, \
\ \ \ \ \ \ \ \ \ \ }c\left( \infty \right) =Q\text{,}  \label{circle}
\end{equation}%
is a circle in the split isotropic Grassmannian \emph{Gr}$_{0}$~$\left(
f\right) $, and all circles have this form.

Moreover, there exist a basis $\mathcal{B}=\left\{ u_{1},\dots 
,u_{n}\right\} $ of $P$ and a map $J:V\rightarrow V$ compactly adapted to $f$
such that  
\begin{equation*}
\gamma \left( s\right) =\text{\emph{span}}~\left\{ \cos s~u_{i}+\sin
s~T\left( u_{i}\right) \mid i=1,\dots ,n\right\}
\end{equation*}
is a geodesic in \emph{Gr}$_{0}$~$\left( f\right) $ satisfying $\gamma 
\left( s\right) =c\left( \tan s\right) $ for all $s\in \mathbb{R}\cup 
\left\{ \infty \right\} $, provided that the metric on \emph{Gr}$_{0}$~$
\left( f\right) $ is invariant by the maximal compact subgroup $\mathcal{U}
^{J}\left( f\right) $ of $\mathcal{U}\left( f\right) $. The relationship 
among $T,\mathcal{B}$ and $J$ is made explicit in the proof. 
\end{theorem}

\noindent \textit{Proof. }The fact that circles have the form (\ref{circle})
is analogous to the similar assertion for the split standard Grassmannians. 
Notice that $c\left( t\right) $ is isotropic if and only if $f\left( 
x+tT\left( x\right) ,y+tT\left( y\right) \right) =0$ for all $x,y\in P$, and
this is equivalent to (\ref{Tskew}), since $P$ and $Q$ are isotropic.

Let $h:P\times P\rightarrow \mathbb{F}$ be defined by $h\left( x,y\right) 
=f\left( x,Ty\right) $. We compute 
\begin{equation*}
h\left( x,y\right) =f\left( x,Ty\right) =-f\left( Tx,y\right) =-\varepsilon
\sigma f\left( y,Tx\right) =-\varepsilon \sigma h\left( y,x\right) \text{,}
\end{equation*}
which is nondegenerate since $T$ is nonsingular.

We consider first the case $\varepsilon =-1$. Then $h$ is bilinear symmetric
or Hermitian symmetric. By the basis theorem \cite{Harvey}, there exists a 
basis $\mathcal{B}=\left\{ u_{1},\dots ,u_{n}\right\} $ of $P$ such that the
matrix of $h$ with respect to $\mathcal{B}$ is diag~$\left( 
I_{k},-I_{n-k}\right) $ for some $k\in \left\{ 0,1,\dots ,n\right\} $ ($
I_{\ell }$ is the $\ell \times \ell $ identity matrix). We call $v_{i}=\rho 
_{i}Tu_{i}$, where $\rho _{i}=1$ for $i=1,\dots ,k$ and $-1$ otherwise. We 
compute 
\begin{equation*}
f\left( u_{i},v_{j}\right) =f\left( u_{i},\rho _{j}Tu_{j}\right) =\rho
_{j}h\left( u_{i},u_{j}\right) =\left( \rho _{j}\right) ^{2}\delta
_{ij}=\delta _{ij}\text{.}
\end{equation*}
Hence $\left\{ v_{1},\dots ,v_{n}\right\} $ is a dual basis of $\left\{ 
u_{1},\dots ,u_{n}\right\} $ with respect to $f$.

Let $J:V\rightarrow V$ be the unique $\bar{\sigma}$-linear map such that $
Ju_{i}=v_{i}$ and $Jv_{i}=\varepsilon u_{i}=-u_{i}$ for all $i=1,\dots ,n$. 
By Lemma \ref{remark} b), $J$ is compactly adapted to $f$ and so $K=_{\text{
def}}\mathcal{U}^{J}\left( f\right) $ is a maximal compact subgroup of $ 
\mathcal{U}\left( f\right) $ by Lemma \ref{Neretin}.

Now, let $S:V\rightarrow V$ be the $\mathbb{F}$-linear map defined by  
\begin{equation}
\left. S\right\vert _{P}=T\text{\ \ \ \ \ \ \ \ \ and\ \ \ \ \ \ \ \ \ }
\left. S\right\vert _{Q}=-T^{-1}\text{.}  \label{SPQ}
\end{equation}
In particular, $Su_{i}=Tu_{i}=\rho _{i}v_{i}$ and $Sv_{i}=-T^{-1}v_{i}=-\rho
_{i}u_{i}$. Let us see that $S\in $ Lie~$\mathcal{U} ^{J}\left( f\right) $.
Indeed,  
\begin{equation*}
JSu_{i}=J\left( \rho _{i}v_{i}\right) =-\rho _{i}u_{i}=Sv_{i}=SJu_{i}\text{.}
\end{equation*}
and one checks similarly that $JSv_{i}=SJv_{i}$ and that $f$ is skew 
symmetric with respect to $f$.

Let $L$ be the isotropy subgroup of $K$ at $P$. Its Lie algebra $\mathfrak{l}
$ consists of the linear maps preserving $P$ and $Q$. Since $S$ interchanges
$P$ and $Q$, tr~$\left( SU\right) =0$ for all $U\in \mathfrak{l}$. This 
means that $S$ belongs to the orthogonal complement of $\mathfrak{l}$ with 
respect to the Killing form. Since $\left( K,L\right) $ is a symmetric pair,
it is well-known that then $\gamma \left( s\right) =\exp \left( sS\right) L$
is a geodesic in $K/L$, provided that the latter carries a $K$-invariant 
Riemannian metric. We have that  
\begin{eqnarray*}
\exp \left( sS\right) \left( u_{i}\right) &=&\cos s~u_{i}+\sin s~T\left(
u_{i}\right) \\
\exp \left( sS\right) \left( T\left( u_{i}\right) \right) &=&-\sin
s~u_{i}+\cos s~T\left( u_{i}\right) \text{.}
\end{eqnarray*}
Therefore, after the canonical identification of $K/L$ with Gr$_{0}\left( 
f\right) ,$ we have that  
\begin{equation*}
\gamma \left( s\right) =\exp \left( sS\right) P=\text{span}~\left\{ \cos
s~u_{i}+\sin s~T\left( u_{i}\right) \mid i=1,\dots ,n\right\} =c\left( \tan
s\right) \text{.}
\end{equation*}

Next we consider the remaining case $\varepsilon =1$. Then $h$ is bilinear 
skew symmetric or Hermitian skew symmetric. By the basis theorem, $n$ is 
even, say $n=2k$, and there exists a basis $\mathcal{B}=\left\{ u_{1},\dots 
,u_{n}\right\} $ of $P$ such that the matrix of $h$ with respect to that 
basis satisfies $h\left( u_{i},u_{k+i}\right) =-h\left( u_{k+i},u_{k}\right)
=1$ for $i=1,\dots ,k$ and the other coefficients are zero. Calling  
\begin{equation*}
v_{i}=T\left( u_{k+i}\right) \text{\ \ \ \ \ \ \ and\ \ \ \ \ \ \ }
v_{k+i}=-T\left( u_{i}\right)
\end{equation*}
for $i=1,\dots ,k$, a straightforward computation yields that $\left\{ 
v_{1},\dots ,v_{n}\right\} $ is a dual basis of $\mathcal{B}$ with respect 
to $f$. Define $S$ as above in (\ref{SPQ}) and $J:V\rightarrow V$ as the 
unique $\bar{\sigma}$-linear map such that $Ju_{i}=v_{i}$ and $
Jv_{i}=u_{i}=\varepsilon u_{i}$ for all $i=1,\dots ,n$. Mutatis mutandis, 
all arguments for the case $\varepsilon =-1$ apply. \hfill $\square $

\subsection{The compact classical Lie groups $SO_{2m}$, $U_{n}$ and $Sp_{n}$}

It is well known that these Lie groups $K$ are, only in another guise, some 
connected components of the split isotropic Grassmannians for a split $ 
\mathbb{F}$-Hermitian form $f$ over $\mathbb{R},\mathbb{C}$ and $\mathbb{H}$
, respectively. In fact, there is bijection between the group $K$ and the 
corresponding Grassmannian:  
\begin{equation*}
\phi :K\rightarrow \text{Gr}_{0}\left( f\right) ,\ \ \ \ \ \ \ \ \ \ \ \phi
\left( A\right) =\text{graph}~\left( A\right) \text{.}
\end{equation*}
One can take $H=\mathcal{U}\left( f\right) $, that is, $O\left( 2m,2m\right)
$, $U\left( n,n\right) $ and $Sp\left( n,n\right) $, acting on $K$ by 
birational transformations: If $\mathcal{A}=\left(  
\begin{array}{cc}
a & c \\ 
b & d%
\end{array}
\right) \in H$, with $a,b,c,d\in \mathbb{F}^{n\times n}$ ($n=2m$ if $\mathbb{%
\ F}=\mathbb{R}$) and $A\in K$, then $\mathcal{A}\cdot A$ is the unique $%
B\in  K $ such that

\begin{equation*}
\left( 
\begin{array}{cc}
a & c \\ 
b & d%
\end{array}
\right) \left( 
\begin{array}{c}
x \\ 
Ax%
\end{array}
\right) \in \text{ graph~}\left( B\right) \text{,}
\end{equation*}
for all $x\in \mathbb{F}^{n}$, that is, $B=\left( bI_{n}+dA\right) \left( 
aI_{n}+cA\right) ^{-1}$.

It easy to check that two operators $S,T\in K$ are opposite if and only if $
Su\neq Tu$ for any unit vector $u\in \mathbb{F}^{2m}$. One may identify $
SO_{2m}$ with the set of positions of an extended (that is, not contained 
in a codimension two subspace) rigid body $B$ in $\mathbb{R}^{2n}$ with the center of 
mass at the origin. Three positions of $B$ turn out to be pairwise opposite 
when passing from one to another no unit vector is fixed. In this case there
is a distinguished circle of positions of $B$ joining them.

We only recall the form of the diagonal geodesics in standard position. They
are $t\mapsto \left( R_{t},\dots ,R_{t}\right) $ for $SO_{2m}$ (where $R_{t}$
is the rotation through angle $t$) and $t\mapsto \left( e^{it},\dots 
,e^{it}\right) $ for the remaining groups.

\subsection{The complex quadric or Grassmannian of oriented planes\label%
{quadric}}

This Grassmannian is very different from the ones in the previous 
subsections. For instance, in contrast to them, two oriented planes which 
are opposite may intersect and even coincide as unoriented planes.

Let $M$ be the Grassmannian of oriented planes in $\mathbb{R}^{2+n}$. We set
$N=2+n$ and consider on $\mathbb{C}^{N}$ the standard bilinear form given by $
\left\langle \left( z_{1},\dots ,z_{N}\right) ,\left( w_{1},\dots 
,w_{N}\right) \right\rangle =z_{1}^{t}w_{1}+\dots +z_{N}^{t}w_{N}$, which is
isomorphic to the one presented in (\ref{table}). Then $M$ can be identified
with the (complex) projectivization $Q$ of the quadric $\left\{ w\in \mathbb{%
\ C}^{N}\mid \left\langle w,w\right\rangle =0\right\} $ as follows:  
\begin{equation}
\psi :M\rightarrow Q\text{,\ \ \ \ \ \ \ \ \ \ \ \ \ \ \ }\psi \left(
u\wedge v\right) =\mathbb{C}\left( u+iv\right) \text{,}  \label{identif}
\end{equation}
for any orthonormal set $\left\{ u,v\right\} $ in $\mathbb{R}^{N}$ (see \cite%
{KN}). The group $SO\left( N,\mathbb{C}\right) $ acts transitively on $Q$, 
in the obvious way. Let $P_{X}$ be the isotropy group of the complex null 
line $\mathbb{C}X$. Its Lie algebra is  
\begin{equation*}
\mathfrak{p}_{X}=\left\{ T\in so\left( N,\mathbb{C}\right) \mid X\text{ is
an eigenvector of }T\right\}
\end{equation*}
and its polar is 
\begin{equation}
\mathfrak{p}_{X}^{\bot }=\left\{ T\in so\left( N,\mathbb{C}\right) \mid TX=0 
\text{ and }T\left( X^{\bot }\right) \subset \mathbb{C}X\right\} \text{,}
\label{pXbot}
\end{equation}
which is abelian, so $\mathfrak{p}_{X}$ has height one. By the list at the 
end of \cite{BP}, $M$ is a self dual symmetric $R$-space.

\smallskip

For the sake of simplicity, in the following we suppose that $N>3$. On the 
one hand, the arguments can be easily adapted to the case $N=3$. On the 
other hand, the action of $SO\left( 3,\mathbb{C}\right) $ above on the 
Grassmannian of oriented planes in $\mathbb{R}^{3}$ is equivalent to the 
well-known action of the direct conformal group on $S^{2}$.

We recall from Lemma 7.124 in \cite{Harvey} the definition of the 
characteristic angles of a pair of oriented planes in $\mathbb{R}^{N}$.

\begin{proposition}
\label{diedral}\emph{\cite{Harvey}}. Given $P,Q$ oriented planes in $\mathbb{%
\ R}^{N}$ with $N>3$, there exist an orthonormal set $\left\{ 
u_{1},v_{1},u_{2},v_{2}\right\} $ in $\mathbb{R}^{N}$ and angles $\alpha 
,\beta $ satisfying 
\begin{equation*}
0\leq \alpha \leq \beta \text{\ \ \ \ \ \ \ \ \ \ \ \ and\ \ \ \ \ \ \ \ \ }
\alpha +\beta \leq \pi
\end{equation*}
$\,$such that $P=u_{1}\wedge u_{2}$ and $Q=\left( \cos \alpha ~u_{1}+\sin 
\alpha ~v_{1}\right) \wedge \left( \cos \beta ~u_{2}+\sin \beta 
~v_{2}\right) $. 
\end{proposition}

The angles $\alpha ,\beta $ are uniquely determined by $P,Q$ and are called 
the characteristic angles of the pair $P,Q$.

The following proposition characterizes the pairs of opposite oriented 
planes and implies that two opposite oriented planes need not be 
complementary; they may even coincide as unoriented subspaces: $u\wedge v$ 
is opposite to $v\wedge u$ for any orthonormal set $\left\{ u,v\right\} $ in
$\mathbb{R}^{N}$.

\begin{proposition}
\label{opposite}Given $P,Q\in M$ with $\psi \left( P\right) =\mathbb{C}X$ 
and $\psi \left( Q\right) =\mathbb{C}Y$ for some non-zero null vectors $
X,Y\in \mathbb{C}^{N}$, the following assertions are equivalent:

\smallskip

\emph{a)} $P$ and $Q$ are opposite.

\smallskip

\emph{b)} $\left\langle X,Y\right\rangle \neq 0$.

\smallskip

\emph{c)} The characteristic angles of the pair $P,Q$ are distinct. 
\end{proposition}

\noindent \textit{Proof. }a) $\Rightarrow $ b) Suppose that $\left\langle 
X,Y\right\rangle =0$. If $X\in \mathbb{C}Y$, then $P=Q$ and so clearly $P$ 
and $Q$ are not opposite. If $X$ and $Y$ are linearly independent, by Lemma  %
\ref{XY} below one can choose an ordered basis $\left\{ X,X^{\prime
},Y,Y^{\prime },Z_{i}\right\} $ such that the associated Gram matrix is
diag~ $\left( R,R,I_{N-4}\right) $, where $R$ is as in (\ref{RR}). Let $T$
be the  linear transformation on $\mathbb{C}^{N}$ defined by $T\left(
X^{\prime }\right) =Y$, $T\left( Y^{\prime }\right) =-X$ and $T$ equal zero
on the  remaining elements of the basis. Then $T\neq 0$ and a
straightforward  computation shows that $T\in \mathfrak{p}_{X}^{\bot }\cap 
\mathfrak{p} _{Y}^{\bot }$. Therefore $P$ and $Q$ are not opposite.

\smallskip

b) $\Rightarrow $ a) Suppose that $\left\langle X,Y\right\rangle \neq 0$. 
Then $X,Y$ together with $X^{\bot }\cap Y^{\bot }$ span $\mathbb{C}^{N}$. If
$T\in \mathfrak{p}_{X}^{\bot }\cap \mathfrak{p}_{Y}^{\bot }$, by (\ref{pXbot}
), $T\left( X\right) =0,$ $T\left( Y\right) =0$ and $T\left( X^{\bot }\cap 
Y^{\bot }\right) \subset \mathbb{C}X\cap \mathbb{C}Y=\left\{ 0\right\} $. 
Hence $T=0$ and so $P$ and $Q$ are opposite.

\smallskip

b) $\Leftrightarrow $ c) Suppose that the characteristic angles of $P$ and $
Q $ are $\alpha ,\beta .$ Then $P$ and $Q$ may be written as in Proposition  %
\ref{diedral} and 
\begin{equation*}
\left\langle u_{1}+iu_{2},\cos \alpha ~u_{1}+\sin \alpha ~v_{1}+i\left( \cos
\beta ~u_{2}+\sin \beta ~v_{2}\right) \right\rangle =\cos \alpha -\cos \beta 
\text{,}
\end{equation*}
from which the assertion follows. \hfill $\square $

\begin{lemma}
\label{XY}Let $X,Y$ be two linearly independent null vectors in $\mathbb{C}
^{N}$ such that $\left\langle X,Y\right\rangle =0$. Then $N>3$ and there 
exist linearly independent null vectors $X^{\prime },Y^{\prime }$ such that 
Gram matrix with respect to the basis $\left\{ X,X^{\prime },Y,Y^{\prime
}\right\} $ is \emph{diag}~$\left( R,R\right) $ with 
\begin{equation}
R=\left( 
\begin{array}{cc}
0 & 1 \\ 
1 & 0%
\end{array}
\right) \text{.}  \label{RR}
\end{equation}
\end{lemma}

\noindent \textit{Proof. }Suppose first that $X^{\bot }\neq Y^{\bot }$.
Since both subspaces have the same dimension $N-1$, there exist vectors $%
U\in X^{\bot }$ and $V\in Y^{\bot }$ with $\left\langle U,Y\right\rangle
=\left\langle V,X\right\rangle =1$. In particular, $X,Y,U,V$ are linearly
independent. Now, a straightforward computation shows that one can choose $%
a,b,c,d\in \mathbb{C}$ such that $X^{\prime }=U+aX+bY$ and $Y^{\prime
}=V+cX+dY$ satisfy the required conditions.

Now suppose that $X^{\bot }=Y^{\bot }.$ Let $X,Y,Z_{3},\dots ,Z_{N}$ be a 
basis of $\mathbb{C}^{N}$, where the first $N-1$ vectors generate $X^{\bot }$, 
and call $A$ the corresponding Gram matrix. Then the first two rows of $A$
have all their coefficients zero, except the last one, in each case. This 
implies that $A$ is singular. A contradiction, since the inner product is 
nondegenerate. \hfill $\square $

\bigskip

Next we describe the circles $c$ in $M$ with $c\left( 0\right) =e_{1}\wedge 
e_{2}$ and $c\left( \infty \right) =e_{2}\wedge e_{1}$. Notice that these 
oriented planes are opposite by Proposition \ref{opposite}, since $
\left\langle e_{1}+ie_{2},e_{2}+ie_{1}\right\rangle =2i\neq 0$; also, their 
characteristic angles are $0,\pi $.

\begin{theorem}
\label{circleQ}Let $P_{0}=e_{1}\wedge e_{2}$ and $Q=e_{2}\wedge e_{1}$. An 
oriented plane $P_{1}$ is opposite to $P$ and $Q$ if and only if it has the 
form 
\begin{equation}
P_{1}=\left( \cos \alpha ~\epsilon_{1}+\sin \alpha ~u\right) \wedge \left(
\cos \beta ~\epsilon_{2}+\sin \beta ~v\right)  \label{P1}
\end{equation}
for some orthonormal sets $\left\{ \epsilon_{1},\epsilon_{2}\right\} $ and $%
\left\{ u,v\right\} $ with  $\epsilon_{1}\wedge\epsilon_{2}=e_{1}\wedge e_{2}
$ and $u,v\in\left\{  e_{1},e_{2}\right\} ^{\bot }$, and some $0\leq \alpha
<\beta $ with $\alpha  +\beta <\pi $. 

The circle $c$ through $P_{0},P_{1}$ and $Q$ is given by $c\left( t\right) 
=\left( u_{t}\wedge v_{t}\right) /\left\Vert u_{t}\right\Vert ^{2}$, with  
\begin{equation*}
u_{t}=\left( 1+t^{2}C\right) ~e_{1}+2ta~u\text{ \ \ \ \ \ \ \ and \ \ \ \ \
\ }v_{t}=\left( 1-t^{2}C\right) ~e_{2}+2tb~v\text{,}
\end{equation*}
where 
\begin{equation}
a=\frac{\sin \alpha }{\cos \alpha +\cos \beta },\ \ \ \ \ \ \ b=\frac{ \sin
\beta }{\cos \alpha +\cos \beta }\ \ \ \ \ \ \ \text{and \ \ \ \ \ \ }C=%
\frac{\cos \alpha -\cos \beta }{\cos \alpha +\cos \beta }\text{.}\ 
\label{ABC}
\end{equation}
\end{theorem}

\medskip

\noindent \textit{Proof. }The first assertion follows from Proposition \ref%
{opposite} and the fact that if $\alpha ,\beta $ are the characteristic 
angles of the pair $P_{0},P_{1}$, then the characteristic angles of $Q,P_{1}$
are $\alpha ,\pi -\beta $. Without loss of generality we may suppose that  $%
\epsilon_{1}=e_{1}$ and $\epsilon_{2}=e_{2}$.

All circles $c$ with $c\left( 0\right) =P_{0}$ and $c\left( \infty \right) 
=Q $ have the form $c\left( t\right) =\exp \left( tZ\right) \cdot P_{0}$ for
some $Z\in \left( \mathfrak{p}_{e_{2}-ie_{1}}\right) ^{\bot }$ with $Z$ 
prevalent. The basis $\mathcal{B}=\left\{ 
e_{1}+ie_{2},e_{1}-ie_{2},e_{3},\dots ,e_{N}\right\} $ has Gram matrix diag~$
\left( 2R,I_{n}\right) $, where $R$ is as in (\ref{RR}). By (\ref{pXbot}), 
the matrix of $Z$ with respect to the basis $\mathcal{B}$ has the form  
\begin{equation}
\left[ Z\right] _{\mathcal{B}}=\left( 
\begin{array}{ccc}
0 & 0 & 0 \\ 
0 & 0 & -z^{t} \\ 
2z & 0 & 0%
\end{array}
\right)  \label{ZB}
\end{equation}
for some $z\in \mathbb{C}^{n}$. A straightforward computation shows that $Z$
is prevalent if and only if $\left\langle z,z\right\rangle \neq 0$. We 
compute  
\begin{equation*}
\left[ \exp tZ\right] _{\mathcal{B}}=\left( 
\begin{array}{ccc}
1 & 0 & 0 \\ 
-t^{2}z^{t}z & 1 & -tz^{t} \\ 
2tz & 0 & 1%
\end{array}
\right) \text{.}
\end{equation*}
Setting $z=xu+iyv$, with $0<x<y$ (hence $\left\langle z,z\right\rangle \neq 
0 $ and so $Z$ is prevalent), we have that 
\begin{eqnarray}
\exp \left( tZ\right) \left( e_{1}+ie_{2}\right) &=&\left(
e_{1}+ie_{2}\right) -t^{2}z^{t}z~\left( e_{1}-ie_{2}\right) +t2z
\label{exptZdot} \\
&=&\left( \left( 1+t^{2}\left( y^{2}-x^{2}\right) \right) e_{1}+t2xu\right)
+i\left( \left( 1-t^{2}\left( y^{2}-x^{2}\right) \right) e_{2}+t2yv\right) 
\text{.}  \notag
\end{eqnarray}

Suppose first that $\alpha >0$. Then, by (\ref{P1}), for $P_{1}=\exp \left( 
Z\right) \cdot \left( e_{1}+ie_{2}\right) $ it suffices that 
\begin{equation}
\cot \alpha =\frac{1+\left( y^{2}-x^{2}\right) }{2x}\ \ \ \ \ \ \ \ \text{
and\ \ \ \ \ \ \ \ \ \ }\cot \beta =\frac{1-\left( y^{2}-x^{2}\right) }{2y} 
\text{.}  \label{cot}
\end{equation}
A straightforward computation shows that $x=a$, $y=b$ is a solution of the 
system of equations (notice that $C=b^{2}-a^{2}$). The case $\alpha =0$ 
follows from (simpler) similar arguments. \hfill $\square $

\begin{proposition}
\label{gamma0}The curve $\gamma _{o}:\mathbb{R}\rightarrow M$ defined by  
\begin{equation*}
\gamma _{o}\left( s\right) =e_{1}\wedge \left( \cos \left( 2\pi s\right)
~e_{2}+\sin \left( 2\pi s\right) ~e_{4}\right) \text{,}
\end{equation*}
is a diametrical geodesic in $M$ endowed with the standard Riemannian 
metric, that is, the one invariant by $SO\left( N\right) $. It satisfies 
that $\gamma _{o}\left( s\right) =c_{o}\left( \tan \left( \pi s\right) 
\right) $, where $c_{o}:\mathbb{R}\cup \left\{ \infty \right\} \rightarrow M$
is the circle through $e_{1}\wedge e_{2}$, $e_{1}\wedge e_{4}$ and $
e_{2}\wedge e_{1}$, that is,  
\begin{equation}
c_{o}\left( t\right) =e_{1}\wedge \left( \frac{1-t^{2}}{1+t^{2}}~e_{2}+\frac{
2t}{1+t^{2}}~e_{4}\right) \text{.}  \label{cSimple}
\end{equation}
\end{proposition}

\noindent \textit{Proof.} Consider the lattices $\Gamma =2\pi \mathbb{Z}^{2}$
and $\Lambda =\pi \left( \mathbb{Z}\left( 1,1\right) +\mathbb{Z}\left( 
-1,1\right) \right) ,$ which contains $\Gamma $. Let $\phi :\mathbb{R}
^{2}\rightarrow M$ be defined by  
\begin{equation*}
\phi \left( s,t\right) =\left( \cos s~e_{1}+\sin s~e_{3}\right) \wedge
\left( \cos t~e_{2}+\sin t~e_{4}\right) \text{,}
\end{equation*}
which can be pushed down to the quotients as indicated in the following 
commutative diagram ($p_{\Gamma }$ and $p$ are the canonical projections). 
The map $\phi _{\Gamma }$ is a double covering of the maximal torus $%
T^{2}=_{ \text{def}}$ Image $\left( \phi _{\Lambda }\right) $ of $M$. 
\begin{equation*}
\begin{array}{ccc}
\mathbb{R}^{2} & \overset{\phi }{\longrightarrow } & M \\ 
\downarrow p_{\Gamma } & \overset{\phi _{\Gamma }}{\nearrow } & \uparrow
\phi _{\Lambda } \\ 
\mathbb{R}^{2}/\Gamma & \overset{p}{\longrightarrow } & \mathbb{R}
^{2}/\Lambda%
\end{array}%
\end{equation*}
Then $t\mapsto \phi \left( 0,t\right) =e_{1}\wedge \left( \cos t~e_{2}+\sin 
t~e_{4}\right) $ is a diagonal geodesic in $T^{2}$, and hence also its 
reparametrization $\gamma _{o}\left( s\right) =\phi \left( 0,2\pi s\right) $
. A straightforward computation shows that $\gamma _{o}\left( s\right) 
=c_{o}\left( \tan \left( \pi s\right) \right) $. By (\ref{exptZdot}), the 
curve $c_{o}$ is the circle through $e_{1}\wedge e_{2}$, $e_{1}\wedge e_{4}$
and $e_{2}\wedge e_{1}$. \hfill $\square $

\bigskip

Now, given any circle $c$ through $P_{0}$ and $Q$, we explicit the 
Riemannian metric on $M$ for which $c$ is a diametral geodesic up to a 
projective reparametrization.

\begin{proposition}
Let $c$ be the circle through $P_{0},P_{1}$ and $Q$ as in Theorem \ref%
{circleQ}. Let $\mathcal{B}$, $a$ and $b$ as in that proposition and suppose
for the sake of simplicity that $u=e_{3}$ and $v=e_{4}$. Let $O\in SO\left(
N, \mathbb{C}\right) $ whose matrix with respect to the basis $\mathcal{B}$
is  
\begin{equation*}
\text{\emph{diag}}~\left( 1/r,r,R_{\sigma },I_{n-2}\right) \text{,}
\end{equation*}
where $a+ib=r\left( \sinh \sigma +i\cosh \sigma \right) $ and  
\begin{equation*}
R_{\sigma }=\left( 
\begin{array}{cc}
\cosh \sigma & -i\sinh \sigma \\ 
i\sinh \sigma & \cosh \sigma%
\end{array}
\right) \text{.}
\end{equation*}
Then $\gamma \left( s\right) =c\left( \tan \left( \pi s\right) \right) $ is 
a diagonal geodesic in a maximal torus of $M,$ provided that $M$ is endowed 
with a Riemannian metric invariant by $O^{-1}~SO_{N}~O$. 
\end{proposition}

\noindent \textit{Proof. }By (\ref{exptZdot}), the curve $c_{o}$ in (\ref%
{cSimple}) coincides up to the identification $\psi $ in (\ref{identif}) 
with $t\mapsto \mathbb{C}\exp \left( tZ\right) \left( e_{1}+ie_{2}\right) $,
where $Z$ is as in (\ref{ZB}) with $z=ie_{4}$. On the other hand, the matrix
of $W=$ Ad~$\left( O\right) ~\left( Z\right) $ with respect to $\mathcal{B}$
has the form (\ref{ZB}) with $z=ae_{3}+ibe_{4}$ and so, again by (\ref%
{exptZdot}), $ c\left( t\right) $ corresponds with $\mathbb{C}\exp \left(
tW\right) \left(  e_{1}+ie_{2}\right) $ via $\psi $. Hence $c\left( t\right)
=O\cdot  c_{0}\left( t\right) $ holds for all $t$. Now, the statement
follows from  Proposition \ref{gamma0}. ~\hfill $\square $

\end{document}